\documentclass[12pt,
amscd]{amsart}
\usepackage{latexsym,amscd,amssymb}  
\theoremstyle{plain}
  \newtheorem{thm}{Theorem}[section]
  \newtheorem{prop}[thm]{Proposition}
  \newtheorem{lem}[thm]{Lemma}
  \newtheorem{cor}[thm]{Corollary}
\theoremstyle{definition}
  \newtheorem{dfn}[thm]{Definition}
  \newtheorem{exmp}[thm]{Example}
\theoremstyle{remark}
  \newtheorem{rem}[thm]{Remark}
\def\ba{{\bf a}}
\def\bb{{\bf b}}

\def\C{C^\bullet}
\def\D{{\mathcal D}^\bullet}
\def\DS{\omega^\bullet}

\def\I{I^\bullet}
\def\cI{{\mathcal I}}

\def\M{M^\bullet}
\def\N{N^\bullet}
\def\P{P^\bullet}
\def\PP{\mathbb P}

\def\rH{\tilde{H}}
\def\cH{\mathcal{H}}
\def\m{{\mathfrak m}}
\def\NN{{\mathbb N}}
\def\const{\underline{k} }
\def\Sq{{\rm Sq}}
\def\ZZ{{\mathbb Z}}
\def\RR{{\mathbb R}}

\def\Hom{\operatorname{Hom}}
\def\Sh{\operatorname{Sh}}
\def\cHom{{\mathcal Hom}}
\def\Or{or}
\def\can{K_{k[Q]}}
\def\Ext{\operatorname{Ext}}
\def\cExt{{\mathcal Ext}}

\def\chara{\operatorname{char}}

\def\ann{\operatorname{Ann}}
\def\supp{\operatorname{supp}}
\def\sp{\operatorname{Sp\acute{e}}}
\def\Supp{\operatorname{Supp}}
\def\<{{\langle}}
\def\>{{\rangle}}
\def\MMZn{{}^*{\rm Mod}}
\def\MMQ{{}^*{\rm Mod}_Q}

\def\too{\longrightarrow}

\def\Id{\operatorname{Id}}
\def\injdim{\operatorname{inj.dim}}
\def\pd{\operatorname{proj.dim}}
\def\relint{\operatorname{rel-int}}
\def\Proj{\operatorname{Proj}}
\numberwithin{equation}{section}
\begin{document}
\title
{Stanley-Reisner rings, sheaves, and Poincar\'e-Verdier duality}
\author{Kohji Yanagawa}
\address{Department of Mathematics, 
Graduate School of Science, Osaka University, Toyonaka, Osaka 
560-0043, Japan}
\email{yanagawa@math.sci.osaka-u.ac.jp}

\maketitle

\begin{abstract}
A few years ago, I defined  a {\it squarefree module} over a polynomial ring 
$S = k[x_1, \ldots, x_n]$ generalizing the Stanley-Reisner ring 
$k[\Delta] = S/I_\Delta$ of a simplicial complex $\Delta \subset 
2^{\{1, \ldots , n\}}$. This notion is very useful in the Stanley-Reisner 
ring theory. In this paper, from a squarefree $S$-module $M$, 
we construct the $k$-sheaf $M^+$ on an $(n-1)$ simplex $B$ which is 
the geometric realization of $2^{\{1, \ldots , n\}}$. 
For example, $k[\Delta]^+$ is (the direct image to $B$ of) the constant sheaf 
on the geometric realization $|\Delta| \subset B$. 
We have $H^i(B, M^+) \cong [H^{i+1}_\m(M)]_0$ for all 
$i \geq 1$. The Poincar\'e-Verdier duality for sheaves $M^+$ on $B$ 
corresponds to the local duality for squarefree modules over $S$. 
For example, if $|\Delta|$ is a manifold, then $k[\Delta]$ is a Buchsbaum 
ring and its canonical module $K_{k[\Delta]}$ is a squarefree module 
which gives the orientation sheaf of $|\Delta|$ with the coefficients in $k$. 
\end{abstract}

\section{Introduction}
This paper presents a new geometric aspect of combinatorial 
commutative algebra on {\it normal semigroup rings.} 
But, in this introduction, we restrict ourselves to 
the polynomial ring case for the simplicity.  In this paper, we use 
the theory of sheaves on a locally compact topological space. 
For this theory, consult \cite{Iver}. 
Basically, we use the same notation as \cite{Iver} here.

Let $S = k[x_1, \ldots, x_n]$ be a polynomial ring over a field $k$, 
and $\Delta$ a simplicial complex whose vertex set is a subset of 
$[n] := \{1, \ldots, n \}$. 
Then the Stanley-Reisner ring $k[\Delta] := S/( \, \prod_{i \in F} x_i 
\mid \text{$F \subset [n]$ with $F \not \in \Delta$} \, )$ of $\Delta$ 
reflects topological properties of the geometric realization $|\Delta|$, 
and has been studied since 1970's (see \cite{BH, St}). 
In \cite{Y}, the author introduced the notion of a {\it squarefree module} 
which generalizes Stanley-Reisner rings. 
This notion allows us to apply homological methods (e.g., derived categories) 
to the theory of Stanley-Reisner rings more systematically. The purpose of 
this paper is to give a geometric meaning of squarefree modules. 

Let $B$ be an $(n-1)$-simplex which is the geometric realization of 
$2^{[n]}$. We construct the $k$-sheaf $M^+$ on $B$ from a squarefree 
module $M$. For example, $k[\Delta]^+ \cong j_* \const_{|\Delta|}$, where 
$\const_{|\Delta|}$ is the constant sheaf on $|\Delta|$ and 
$j: |\Delta| \to B$ is the embedding map. Let $\Sq$ be the category of 
squarefree modules, and $\Sh(B)$ the category of $k$-sheaves on $B$.   
Then the functor $(-)^+: \Sq \to \Sh(B)$ is exact. 

If $M$ is a squarefree $S$-module, Theorem~\ref{cohomology} gives 
an isomorphism 
$$H^i(B, M^+) \cong [H_\m^{i+1}(M)]_0 \quad  \text{for all $i \geq 1$},$$ 
and an exact sequence 
$$0 \to [H_\m^0(M)]_0 \to M_0 \to H^0( B, M^+) \to [H_\m^1(M)]_0 \to 0,$$
where $H_\m^i(-)$ stands for the local cohomology module with support in 
the maximal ideal $\m := (x_1, \ldots, x_n)$. So our functor $(-)^+$ is 
somewhat analogous to ``$\Proj$" of the scheme theory. 

If $\ba = (a_1, \ldots, a_n) \not \in -\NN^n$, 
it is well-known that $[H_\m^i(M)]_{\ba} = 0$. 
If $0 \ne \ba \in -\NN^n$, $[H_\m^i(M)]_{\ba}$ is isomorphic to the cohomology 
with compact support $H^{i-1}_c(U_F, j^*M^+)$, where $U_F$ is the open subset 
determined by $F := \{ i \in [n] \mid a_i < 0 \}$ 
(see Theorem~\ref{cohomology2} for detail) 
and $j: U_F \to B$ is the embedding map.  These results generalize 
a well-known formula of Hochster on $H_\m^i(k[\Delta])$. 

Let $D^b(\Sq)$ be the bounded derived category of $\Sq$ and 
$\DS_S \in D^b(\Sq)$ an injective resolution of $K_S[n-1]$. Set 
$\DS_{k[\Delta]} := \Hom_{k[\Delta]}^\bullet(k[\Delta], \DS_S)$. 
Then $\DS_{k[\Delta]}$ is a complex of squarefree $k[\Delta]$-modules, and
isomorphic to a (non-normalized) $\ZZ^n$-graded dualizing complex of 
$k[\Delta]$ in the derived category. 
Let $\D_{|\Delta|}$ be a dualizing complex 
of the topological space $|\Delta|$ with the coefficients 
in $k$. Corollary~\ref{dualizings} states that 
$\D_{|\Delta|} \cong j^*(\DS_{k[\Delta]})^+$ in  $D^b(\Sh(|\Delta|))$, 
where $j:|\Delta| \to B$ is the embedding map. 
(Since the functors $(-)^+: \Sq \to \Sh(B)$ and 
$j^*: \Sh(B) \to \Sh(|\Delta|)$ are exact, we have the functor  
$j^*(-)^+: D^b(\Sq) \to D^b(\Sh(|\Delta|))$.)
Moreover, if each component $M^i$ of $\M \in D^b(\Sq)$ is a 
$k[\Delta]$-module, we have  
$$R \cHom_{\Sh(|\Delta|)}(\, j^*(\M)^+, \D_{|\Delta|} \,) \cong
j^* (R\Hom_{k[\Delta]} ( \, \M, \DS_{k[\Delta]})^+)$$
in $D^b(\Sh(|\Delta|))$. See Theorem~\ref{local Verdier duality}. 
If further $[H^i(\M)]_0 = 0$ for all $i$ (note that the sheaf $M^+$ 
does not reflect the degree 0 component $M_0$), we have 
$$\Ext^i_{\Sh(|\Delta|)}(\, j^*(\M)^+, \, \D_{|\Delta|} \, )  
\cong [\Ext^i_{k[\Delta]}(\, \M, \, \DS_{k[\Delta]} \, )]_0.$$
So the Poincar\'e-Verdier duality for $|\Delta|$ corresponds the 
local duality for $k[\Delta]$ in our context. 
For example, if $|\Delta|$ is a manifold, then $k[\Delta]$ is a Buchsbaum 
ring and the canonical module $K_{k[\Delta]}$ of $k[\Delta]$ is a squarefree 
module which gives the orientation sheaf of $|\Delta|$ with the coefficients 
in $k$. So, the well-known duality between $H_\m^i(k[\Delta])$ and 
$H_\m^j(K_{k[\Delta]})$ (\cite[II. Theorem~4.9]{SV})
corresponds to the Poincar\'e duality for $|\Delta|$. 

\section*{Acknowledgments}
I am grateful to Professor Ezra Miller. This paper 
grew out of the discussion with him about his paper \cite{Mil}.  
I also thank Professor Takesi Kawasaki for giving me useful comments 
on (graded) Buchsbaum modules.

\section{Preliminaries}
Let $Q \subset \NN^n$ be an affine semigroup (i.e., a finitely generated 
sub-semigroup containing 0), and $k[Q] = k[ \, x^\ba \mid \ba \in Q \, ] 
\subset S := k[\, x_1, \ldots, x_n \, ]$ 
the semigroup ring of $Q$ over a field $k$. 
Here $x^\ba$ for $\ba = (a_1, \ldots, a_n) \in \NN^n$ 
means the monomial $\prod x_i^{a_i} \in S$. We always assume that 
$Q$ is saturated (i.e., if $\ba \in \NN^n$ satisfies $m\ba \in Q$ for a 
positive integer $m$, then $\ba \in Q$) and $\ZZ Q = \ZZ^n$. Thus $k[Q]$ is a 
normal Cohen-Macaulay $\ZZ^n$-graded ring of dimension $n$ with 
the graded maximal ideal $\m = ( \, x^\ba \mid 0 \ne \ba \in Q \, )$. For 
basic properties of ($\ZZ^n$-graded modules over) $k[Q]$ and the 
related notions from convex geometry, see \cite{BH, GW}. 
 
Consider the polyhedral cone  $\RR_{\geq 0}Q \subset \RR^n$ spanned by 
$Q \, (\subset \NN^n \subset \RR^n)$. 
Let $L$ be the set of non-empty faces of $\RR_{\geq 0}Q$. 
The order by inclusion makes $L$ a finite poset. 
If $p \in \RR_{\geq 0} Q$, there is a unique face $F \in L$ 
such that the relative interior $\relint(F)$ of $F$ contains $p$. 
We call this $F$ the {\it support} of $p$, and denote it by $\supp(p)$. 

If $\RR_{\geq 0}Q$ is spanned by $n$ vectors as a polyhedral 
cone, we say $k[Q]$ is {\it simplicial}. In this case, $L$ is isomorphic 
to the boolean lattice $2^{[n]}$ as a poset. 
For example, the polynomial ring $k[\NN^n] = k[x_1, \ldots, x_n]$ is 
a simplicial semigroup ring. 

Let $H$ be a hyperplane of $\RR^n$ which intersects the cone $\RR_{\geq 0}Q$ 
transversally. Consider the $(n-1)$-dimensional polytope 
$B:= H \cap \RR_{\geq 0}Q$. 
Of course, $B$ is homeomorphic to a closed ball of dimension $n-1$.  
If $k[Q]$ is simplicial, then $B$ is a simplex. For 
a face $F \in L$, set $|F| := F \cap H$ to be the face of $B$, 
and $|F|^\circ := \relint(|F|)$ its relative interior. 
If $\Delta \subset L$ is an {\it order ideal} (i.e., $F \in \Delta$, 
$G \in L$ and $F \supset G$ $\Rightarrow$ $G \in \Delta$), then 
$|\Delta|  := \coprod_{F \in \Delta} |F|^\circ$ is a finite regular cell 
complex. 

For a $\ZZ^n$-graded $k[Q]$-module $M$ and $\ba \in \ZZ^n$, $M_\ba$ 
denotes the degree $\ba$ component of $M$. 
Let $\MMZn$ be the category of $\ZZ^n$-graded $k[Q]$-modules. 
Here a morphism $f$ in $\MMZn$ is a $k[Q]$-homomorphism 
$f: M \to N$ with $f(M_\ba) \subset N_\ba$ for all $\ba \in \ZZ^n$.  

We assign an order ideal $\Delta \subset L$ to the  
ideal $I_\Delta := (\, x^\ba \mid \text{$\ba \in Q$ and 
$\supp(\ba) \not \in \Delta$} \, )$ of $k[Q]$. 
Set $k[\Delta] := k[Q]/I_\Delta$. Clearly, 
$$
k[\Delta]_\ba \cong 
\begin{cases}
k & \text{if $\ba \in Q$ and $\supp(\ba) \in \Delta$,}\\
0 & \text{otherwise.}
\end{cases}
$$
In particular, if $\Delta = L$ (resp. $\Delta = \emptyset$), 
then $I_\Delta = 0$ (resp. $I_\Delta = k[Q]$) and $k[\Delta] = k[Q]$ 
(resp. $k[\Delta]=0$). If $\Delta \ne \emptyset$ or $\{\, \{0\} \, \}$, 
then $\dim k[\Delta] = \dim |\Delta| +1$, 
where $\dim |\Delta|$ is the dimension as a cell complex. 
When $k[Q]$ is a polynomial ring, $k[\Delta]$ is nothing other than 
the Stanley-Reisner ring of a simplicial complex $\Delta$. 
(If $k[Q]$ is simplicial, $\Delta$ can be seen as a simplicial complex, 
and $|\Delta| = \coprod_{F \in \Delta} |F|^\circ$ is homeomorphic to 
the geometric realization of $\Delta$ as a simplicial complex.) 

We now recall the definition of {\it squarefree} $k[\Delta]$-modules. 

\begin{dfn}[\cite{Y, Y2}]
A $\ZZ^n$-graded $k[Q]$-module $M= \bigoplus_{\ba \in \ZZ^n} M_\ba$ 
is {\it squarefree}, if the following two conditions are satisfied. 
\begin{itemize}
\item[(1)] $M$ is finitely generated and 
$Q$-graded (i.e.,  $M_{\ba} = 0$ for all $\ba \not \in Q$).  
\item[(2)] The multiplication map $M_\ba \ni y \mapsto x^\bb y 
\in M_{\ba + \bb}$ is bijective for all 
$\ba, \bb \in Q$ with $\supp (\ba+\bb) = \supp (\ba)$. 
\end{itemize}
\end{dfn} 

The $\ZZ^n$-graded canonical module $\can$ of $k[Q]$ is a squarefree module. 
In fact, $\can$ is isomorphic to the ideal 
$(\, x^\ba \mid \text{$\ba \in Q$ with $\supp(\ba) = \RR_{\geq 0} Q$} \, )$ 
of $k[Q]$. The quotient rings $k[\Delta]$ (in particular, $k[Q]$ itself) 
are also squarefree. 

If $M$ is squarefree, 
then $M_\ba \cong M_\bb$ for all $\ba, \bb \in Q$ with $\supp(\ba) 
= \supp(\bb)$. In fact, since $\supp(\ba)= \supp (\ba +\bb) 
= \supp(\bb)$, we have $M_\ba \cong M_{\ba + \bb} \cong M_\bb$.

Denote the full subcategory of $\MMZn$ consisting of all squarefree 
$k[Q]$-modules by $\Sq$. For $M \in \MMZn$ and $\ba \in \ZZ^n$, 
$M(\ba)$ denotes the shifted module of $M$ with $M(\ba)_\bb = M_{\ba+\bb}$. 
If $M, N \in \MMZn$ and $M$ is finitely generated, then $\Hom_{k[Q]}(M,N)$ 
has the natural $\ZZ^n$-grading with $[\Hom_{k[Q]}(M,N)]_\ba = 
\Hom_{\MMZn}(M,N(\ba))$. 

\begin{lem}[{\cite[\S4]{Y2}}]\label{basic} 
\begin{itemize}
\item[(1)] $\Sq$ is a thick abelian subcategory of $\MMZn$ (i.e., closed under 
kernels, cokernels, and extensions in $\MMZn$). 
\item[(2)] $\Sq$ is an abelian category with 
enough projectives and injectives. An indecomposable 
projective (resp. injective) object 
in $\Sq$ is isomorphic to $$J_F : = (x^\ba \mid \text{$\ba \in Q$ with 
$\supp(\ba) \supset F$}) \subset k[Q]$$ 
$$\text{( \, resp. 
$k[F] := k[Q]/(x^\ba \mid \text{$\ba \in Q$ with 
$\supp(\ba) \not \subset F$})$ \, )} $$ 
for some $F \in L$. And both $\pd_{\Sq} M$ and 
$\injdim_{\Sq} M$ are at most $n$ for all $M \in \Sq$. 
\item[(3)] The projective object 
$J_F$ is a Cohen-Macaulay $k[Q]$-module of dimension $n$. And  
$$\Hom_{k[Q]}(J_F, \can) \cong 
(\, x^\ba \mid \text{$\ba \in Q$ such that 
$\supp(\ba) \vee F = \RR_{\geq 0} Q$} \, ) ,$$  
where $\supp(\ba) \vee F \in L$ is the  
smallest face containing both $\supp(\ba)$ and $F$. 
In particular, $\Hom_{k[Q]}(J_F, \can)$ is squarefree again. 
\end{itemize}
\end{lem}

For derived categories, we use the same notation as \cite{RD} 
(unless otherwise specified). 
In particular, for a module $M$ and an integer $i$, $M[i]$ means 
the complex $\cdots \to 0 \to M \to 0 \to \cdots $ 
with $M$ at the $(-i)^{\rm th}$ place. 

\begin{lem}\label{full subcat}
We have the canonical category equivalence $D^b(\Sq) \cong D^b_\Sq(\MMZn)$, 
and $D^b(\Sq)$ can be seen as a full subcategory of $D^b(\MMZn)$. 
\end{lem}

\begin{proof}
Let $\MMQ$ be the full subcategory of $\MMZn$ consisting of finitely generated 
$Q$-graded modules. Then $\MMQ$ is a thick abelian subcategory of $\MMZn$. 
Moreover,  
$\MMQ$ has enough projectives, and projective objects $k[Q](-\ba)$ with  
$\ba \in Q$ are also projective in $\MMZn$. Thus $D^b_{\MMQ}(\MMZn) \cong 
D^b(\MMQ)$ and $D^b_{\MMQ}(\MMZn)$ is a full subcategory of $D^b(\MMZn)$. 
On the other hand, $\Sq$ is a thick abelian subcategory of $\MMQ$, and 
an injective object $k[F]$ of $\Sq$ is also injective in $\MMQ$ by 
\cite[Remark~2.5]{Mil}. 
So $D^b(\Sq) \cong D^b_{\Sq}(\MMQ)$, and $D^b_{\Sq}(\MMQ)$ is a full 
subcategory of $D^b(\MMQ)$, which can be seen as a full subcategory of 
$D^b(\MMZn)$. 
\end{proof}

Next we will study $R\Hom_{k[Q]}(\M, \can)$ for a complex $\M \in D^b(\MMZn)$. 
Here $R\Hom_{k[Q]}(\M, \can)$ is the ``$R\Hom$" as complexes of (non-graded) 
$k[Q]$-modules. But if each $H^i(\M)$ is finitely generated, 
$R\Hom_{k[Q]}(\M, \can)$ has a natural $\ZZ^n$-grading, and 
defines an object in $D^b(\MMZn)$. 
In fact, if $\P$ is a $\ZZ^n$-graded finite free resolution  
of $\M$, then $R\Hom_{k[Q]}(\M, \can) \cong \Hom_{k[Q]}^\bullet(\P, \can)$ and 
each $\Hom_{k[Q]}^i(\P, \can) \, (\, = \Hom_{k[Q]}(P^{-i}, \can) \,)$ 
has the $\ZZ^n$-grading. We can also define $R\Hom_{k[Q]}(\M, \can)$ using 
a $\ZZ^n$-graded injective resolution $\I \in D^b(\MMZn)$ of $\can$, 
but we get the same $\ZZ^n$-grading. 

\begin{lem}\label{local dual}
If $\M \in D^b_{\Sq}(\MMZn)$, then $R\Hom_{k[Q]}(\M, \can)$ is in 
 $D^b_{\Sq}(\MMZn)$ too. That is, $\Ext^i_{k[Q]}(\M, \can)$ 
 is squarefree for all $i$. 
\end{lem}

\begin{proof}
By Lemma~\ref{full subcat}, 
we may assume that $\M \in D^b(\Sq)$. Then we have a projective 
resolution $\P \in D^b(\Sq)$ of $\M$. By Lemma~\ref{basic} (3),  
$\Ext^i_{k[Q]}(P^j, \can) = 0$ for all $i \ne 0$ and all $j$.  Hence 
we have $R\Hom_{k[Q]}(\M, \can) \cong \Hom_{k[Q]}^\bullet(\P,\can)$. 
But each $\Hom_{k[Q]}(P^j,\can)$ is squarefree by Lemma~\ref{basic} (3). 
\end{proof}

Take some $\ba(F) \in Q \cap \relint(F)$ for each $F \in L$. 
For a squarefree module $M$, set $M_F := M_{\ba(F)}$. 
If $F, G \in L$ and $G \supset F$, \cite[Theorem~3.3]{Y2} gives a 
$k$-linear map $\varphi^M_{G, F}: M_F \to M_G$.  These maps satisfy 
$\varphi^M_{F,F} = \Id$ and 
$\varphi^M_{G, F} \circ \varphi^M_{F, E} = \varphi^M_{G,E}$ for all 
$G \supset F \supset E$. For $F \in L$, we define the complex 
$\C_F(M) : 0 \to C_F^0 \to C_F^1 \to \cdots \to 
C_F^n \to 0$ of $k$-vector spaces by 
$$C_F^i = \bigoplus_{\substack{G \in L, \, G \supset F \\ \dim G = i}} M_G,$$
and the differential 
$$d : C_F^{i} \supset M_G \ni y 
\longmapsto \sum_{\substack{G' \in L, \, G' \supset G\\ 
\dim G' = i +1}} \varepsilon(G', G) \cdot \varphi^M_{G',G} (y) \ \in 
\bigoplus_{\substack{G' \in L, \, G' \supset G\\ \dim G' = i +1}} 
M_{G'} \subset C_F^{i+1}.$$
Here $\varepsilon$ is an incidence function on the cell complex 
$B = \coprod_{F \in L} |F|^\circ$. 
The complex $\C_F(M)$ does not depend on the particular choice of $\ba(F)$'s 
up to isomorphism. 
By the computation of a \v{C}ech complex with supports in $\m$, 
we have the following.

\begin{lem}[{\cite[Theorem~3.10]{Y2}}]\label{C_F} 
Let the notation be as above. If $\ba \not \in Q$, then 
$[H_\m^i(M)]_{-\ba} = 0$. If $\ba \in Q$ and $\supp(\ba) = F$, 
then $[H_\m^i(M)]_{-\ba} \cong H^i(\C_F(M))$. 
\end{lem}

\section{Sheaves associated with squarefree modules}
We keep the same notation as above. For a squarefree module $M$, 
set $$\sp(M) := \coprod_{F \in L} |F|^\circ \times M_F.$$ 
Let $\pi : \sp(M) \to B$ be the projection map which sends $(p, m) \in 
|F|^\circ \times M_F \subset \sp(M)$ to $p \in |F|^\circ \subset B$. 
For an open subset $U \subset B$ and a map $s: U \to \sp(M)$, 
we will consider the following conditions:  

\begin{itemize}
\item[$(*)$]  $\pi \circ s = \Id_{U}$ and $s_q = \varphi^M_{G, F}(s_p)$
for all $p, q \in U$ such that $F := \supp (p)$ is contained in  
$G := \supp(q)$. Here $s_p$ (resp. $s_q$) is the element of $M_F$ 
(resp. $M_G$) with $s(p) = (p, s_p)$ (resp. $s(q) = (q, s_q)$).  
\item[$(**)$] There is an open covering $U = \bigcup_{\lambda \in \Lambda} 
U_\lambda$ such that the restriction of $s$ to $U_\lambda$ satisfies $(*)$ for 
all $\lambda \in \Lambda$. 
\end{itemize}

Now we define the $k$-sheaf associated to $M$ on $B$, denoted by $M^+$, 
as follows. The sections $M^+(U)$ of $M^+$ over an open set $U$ is 
$$\{ \, s \mid \text{$s: U \to \sp(M)$ is a map satisfying $(**)$} \,\}$$
and the restriction map $M^+(U) \to M^+(V)$ is the natural one. 
(That $M^+$ is actually a sheaf is obvious.)

We say an open set $U$ of $B$ is {\it neat} with respect to a face $F \in L$, 
if $U$ itself and $U \cap |G|^\circ$ are connected for all $G \in L$ 
with $G \supset F$, and $q \in U$ implies $\supp(q) \supset F$. 
For example, the open set $U_F := \coprod_{\, G \supset F} 
|G|^\circ$ is neat with respect to $F$. For $x \in |F|^\circ$ and 
sufficiently small $\varepsilon > 0$, 
$U_{\varepsilon}(x) := \{ \, y \in B \mid d(x, y) < \varepsilon \, \}$  
is also neat with respect to $F$, 
where $d(-,-)$ stands for the usual metric of $\RR^n \, (\supset B)$. 
We can easily check the following. 
\begin{itemize}
\item[(i)] Assume that $U \cap |F|^\circ$ is connected, and 
let $s \in M^+(U)$ be a section. Then there is some $y \in M_F$ such that 
$s(p) = (p, y)$ for all $p \in U \cap |F|^\circ$. 
\item[(ii)] Assume that $U$ is neat with respect to $F$.  For any $y \in M_F$, 
the map $s_y : U \to \sp(M)$ defined by 
$(U \cap |G|^\circ) \ni p \mapsto (p, \varphi_{G,F}(y))$ satisfies $(*)$. 
In particular, $s_y \in M^+(U)$. 
\item[(iii)] If $U$ is neat with respect to $F$, 
any section $s \in M^+(U)$ coincides with $s_y$ of (ii) for some 
$y \in M_F$. 
\end{itemize}
Hence, if $U$ is neat with respect to $F$, then $M^+(U) \cong M_F$. 
Note that the set of neat open sets is an open base of $B$. 
Thus, for a point $p \in |F|^\circ$, the stalk $(M^+)_p$ of $M^+$ at 
$p$ is isomorphic to $M_F$. So $\sp(M)$ is the etale space of the sheaf 
$M^+$. 

\medskip

Let $\Psi \subset L$ be an order filter of the poset $L$, that is, 
$F \in \Psi$, $G \in L$, and $G \supset F$ imply $G \in \Psi$. Then 
$U_\Psi := \coprod_{F \in \Psi} |F|^\circ$ is an open subset of $B$. 
If $M$ is a squarefree module, then the submodule 
$$M_\Psi := \bigoplus_{\ba \in Q, \, \supp_+ (\ba) \in \Psi} M_\ba$$ 
is also squarefree. Moreover, we have the following. 

\begin{lem}\label{truncation} 
The sheaf $(M_\Psi)^+$ is isomorphic to $j_! j^* M^+$, where $j: U_\Psi \to B$ 
is the embedding map. 
\end{lem} 

\begin{proof}
Straightforward. 
\end{proof}

\begin{exmp}\label{constant}
(1) Let $\Delta \subset L$ be an order ideal, and $j$ the embedding map 
from the closed subset $|\Delta| = \coprod_{F \in \Delta} |F|$ to $B$. 
Then the  sheaf $k[\Delta]^+$ is isomorphic to $j_* \const_{|\Delta|}$, where 
$\const_{|\Delta|}$ is the constant sheaf on $|\Delta|$. 

(2) Let $J_F$ be the projective object in $\Sq$ associated with a face 
$F \in L$. Then the sheaf $(J_F)^+$ is isomorphic to 
$j_! \const_{U_F}$, where $j$ is the embedding map from the open set 
$U_F = \coprod_{G \in L, \, G \supset F} |G|^\circ$ to $B$. Note that 
$$U_F \cong 
\begin{cases}
\RR^{n-1}, & \text{if $F = \RR_{\geq 0} Q$,} \\
\RR_+^{n-1} := \{ \, (y_1, \ldots, y_{n-1}) \in \RR^{n-1} \mid 
y_{n-1} \geq 0 \, \}, & \text{if $F \ne \RR_{\geq 0} Q, \{ 0 \}$,} \\
B^{n-1} := \{ \, (y_1, \ldots, y_{n-1}) \in \RR^{n-1} \mid 
\sum_{i=1}^{n-1} y_i^2 \leq 1 \, \}, & \text{if $F = \{ 0 \}$.}
\end{cases}$$

(3) Let $\Delta, \Sigma \subset L$ be order ideals with 
$\Delta \supset \Sigma$. We have $I_\Delta \subset I_\Sigma$. 
Set $I_{\Delta / \Sigma} := I_\Sigma/I_\Delta$ (see \cite[III.7]{St}). 
If $\Sigma = \emptyset$ (resp. $\Delta = L$), then 
$I_{\Delta / \Sigma} = k[\Delta]$ (resp. $I_{\Delta / \Sigma} = I_\Sigma$). 
It is easy to see that $I_{\Delta / \Sigma}$ is a squarefree module with 
$(I_{\Delta / \Sigma})^+ \cong j_! \, \const_{|\Delta| - |\Sigma|}$, 
where $j$ is the embedding map from the locally closed subset 
$|\Delta| - |\Sigma|$ to $B$.
\end{exmp}

 For a topological space $X$, $\Sh(X)$ denotes the category of $k$-sheaves 
on $X$ (i.e., the category of $\const_X$-modules).

If $M$ is a squarefree module, $M_{> 0}$ denotes the submodule 
$\bigoplus_{\ba \in Q \setminus\{0\}} M_\ba$ of $M$. 
Then $M_{>0}$ is squarefree again, 
and $M^+ \cong (M_{>0})^+$. For a complex $0 \to L \to M \to N \to 0$ 
of squarefree modules, the complex of sheaves 
$0 \to L^+ \to M^+ \to N^+ \to 0$ is exact if and only if 
$0 \to L_F \to M_F \to N_F \to 0$ is exact for all $\{ 0\} \ne F 
\in L$. Hence the functor $(-)^+: \Sq \to \Sh(B)$ is exact. 
But this functor is neither full nor faithful. 
The degree 0 component $M_0$ causes this problem. 
Let $\Sq_+$ be the full subcategory of $\Sq$ consisting of all $M$ with 
$M_0 = 0$. It is easy to see that 
the functor $(-)^+: \Sq_+ \to \Sh(B)$ is fully faithful. 

\begin{thm}\label{cohomology}
If $M$ is a squarefree $k[Q]$-module, we have an isomorphism 
$$H^i(B, M^+) \cong [H_\m^{i+1}(M)]_0 \quad  \text{for all $i \geq 1$},$$ 
and an exact sequence 
\begin{equation}\label{lower dimenion}
0 \to [H_\m^0(M)]_0 \to M_0 \to H^0( B, M^+) \to [H_\m^1(M)]_0 \to 0.
\end{equation}
In particular, if $M \in \Sq_+$, then 
$H^i(B, M^+) \cong [H_\m^{i+1}(M)]_0$ for all $i \geq 0$. 
\end{thm}

\begin{proof}
As usual, let $\Gamma_\m: \MMZn \to \MMZn$ be the functor defined by 
$\Gamma_\m(N) := \{ y \in N \mid \text{$\m^l y = 0$ for $l \gg 0$}\}$, 
and $\Gamma(B, -): \Sh(B) \to \operatorname{vect}_k$ the global sections 
functor. 

Let $I^\bullet$ (resp.  $\check{I}^\bullet$)
be a minimal injective resolution of $M$ in $\Sq$ (resp. in $\MMZn$), 
and consider the exact sequence 
\begin{equation}\label{complexes}
0 \to \Gamma_\m (I^\bullet) \to I^\bullet \to I^\bullet/ \Gamma_\m(I^\bullet) 
\to 0
\end{equation}
of cochain complexes. Put $J^\bullet := I^\bullet/ \Gamma_\m(I^\bullet)$. 
Each component of $J^\bullet$ is a direct sum of copies of 
$k[F]$ for various $\{ 0 \} \ne F \in L$. 
Since $k[F]^+$ is the constant sheaf on $|F|$ which is homeomorphic to 
a closed ball, we have $H^i(B, k[F]^+) = H^i(|F|; k)= 0$ for all $i \geq 1$.  
Hence $(J^\bullet)^+ \, (\, \cong (I^\bullet)^+ \,)$ 
gives a $\Gamma(B, -)$-acyclic resolution of $M^+$. 
It is easy to see that $[J^\bullet]_0 \cong \Gamma(B, (J^\bullet)^+)$. 
By \cite[Theorem~2.4]{Mil}, $\I$ coincides with the $Q$-graded part 
$\bigoplus_{\ba \in Q} [\check{I}^\bullet]_\ba$ of $\check{I}^\bullet$. 
Thus we have $[H^i(\Gamma_\m(I^\bullet))]_0 = 
[H^i(\Gamma_\m(\check{I}^\bullet))]_0 = [H_\m^i(M)]_0$. 
So the first and the second assertions follow from \eqref{complexes}, since
$[H^0(I^\bullet)]_0 \cong M_0$ and $H^i(I^\bullet) = 0$ for all $i \geq 1$. 

To prove the last isomorphism, we may assume that $i = 0$. 
But the isomorphism follows from the exact sequence \eqref{lower dimenion}, 
since $H_{\m}^0(M) = M_0 = 0$ in this case. 
\end{proof}

\begin{rem}  
Let $M$ be a finitely generated $\ZZ$-graded module over 
$S= k[x_1, \ldots, x_n]$. 
Then we have an algebraic coherent sheaf $\tilde{M}$ on $\PP^{n-1} = 
\Proj (S)$. Like our functor $(-)^+$, if $\dim_k M < \infty$, then 
$\tilde{M} = 0$. Moreover, it is well-known that 
$H^i(\PP^{n-1}, \tilde{M}) \cong [H_\m^{i+1}(M)]_0$ 
for all $i \geq 1$,  and $$0 \to [H_\m^0(M)]_0 \to M_0 \to 
H^0(\PP^{n-1} , \tilde{M}) \to [H_\m^1(M)]_0 \to 0 \qquad \text{(exact),}$$ 
(see, for example, \cite[p.38]{SV}). 
So Theorem~\ref{cohomology} gives an analogy between $\Proj$ and 
our $(-)^+$.  
\end{rem}

Recall that we chose $\ba(F) \in Q \cap F$ for each $F \in L$ 
in the previous section. 
By the graded local duality, the $\ZZ^n$-graded $k$-dual of $H_\m^i(M)$ 
is isomorphic to the squarefree module $\Ext^{n-i}_{k[Q]}(M, \can)$. 
So to determine the $\ZZ^n$-graded Hilbert function of $H_\m^i(M)$, 
it suffices to know $[H_\m^i(M)]_{-\ba(F)}$ for each $F$. Since 
Theorem~\ref{cohomology} deal with the case when $F = \{ 0 \}$ 
(i.e., $\ba(F) = 0$), we may assume that $F \ne \{ 0 \}$.

\begin{thm}\label{cohomology2}
Let $M$ be a squarefree $k[Q]$-module, and $j$ the embedding map from 
the open set $U_F = \coprod_{G \supset F} |G|^\circ$ to $B$.  
If $F \ne \{ 0\}$, we have 
$$H^i_c(U_F, j^*M^+) \cong [H_\m^{i+1}(M)]_{-\ba(F)} 
\quad \text{for all $i \geq 0$},$$ 
where $H_c^i(-)$ stands for the cohomology with the compact support. 
\end{thm}

\begin{proof}
Let $\Psi := \{ G \in L \mid G \supset F \}$ be the order filter of $L$. 
Under the same notation as Lemma~\ref{truncation}, we have $U_\psi = U_F$. 
We have the following. 
\begin{eqnarray*}
[H_\m^{i+1}(M)]_{-\ba(F)} 
&\cong& H^{i+1}(C_F^\bullet(M)) \qquad \text{(by Lemma~\ref{C_F})}\\
&\cong& H^{i+1}(C^\bullet_{ \{ 0 \} } (M_\Psi)) \\
&\cong& [H^{i+1}_\m(M_\Psi)]_0 \qquad \text{(by Lemma~\ref{C_F})}\\
&\cong& H^i(\, B, \, (M_\Psi)^+ \, ) \qquad 
\text{(by Theorem~\ref{cohomology}. Note that $M_\Psi \in \Sq_+$)}\\ 
&\cong& H^i(\, B, \, j_! j^* M^+\, ) \qquad 
\text{(by Lemma~\ref{truncation})}\\                      
&\cong& H^i_c(\, U_F, \, j^* M^+\, ) \qquad 
\text{(by \cite[III, Corollary~7.3]{Iver})}\\ 
\end{eqnarray*}
\end{proof}

\begin{rem}\label{Hoch}
When $k[Q]$ is a polynomial ring $k[\NN^n] = k[x_1, \ldots, x_n]$,  
Theorems~\ref{cohomology} and \ref{cohomology2} generalize a well-known 
formula of Hochster (\cite[Theorem~5.3.8]{BH}, see also 
\cite[Lemma~5.4.5]{BH}). This formula states that 
$[H_\m^{i+1}(k[\Delta])]_0 \cong \rH_i(|\Delta| ; k)$ for all $i \geq 0$, 
where the right hand side is the $i$th reduced homology group of $|\Delta|$. 
On the other hand, 
Theorem~\ref{cohomology} states that $[H_\m^{i+1}(k[\Delta])]_0 = 
H^i(B, k[\Delta]^+)$ for all $i \geq 1$ and $H^0(B, k[\Delta]^+) 
= [ H_\m^1(k[\Delta]) ]_0 \oplus k[\Delta]_0 \cong 
[  H_\m^1(k[\Delta]) ]_0 \oplus k.$ But  
$H^i(B, k[\Delta]^+) = H^i( B, \, j_* \const_{|\Delta|} ) 
= H^i( |\Delta| ; \, k  ) = \rH_i( |\Delta|; \, k  )$  
for all $i \geq 1$, and $H^0(B, k[\Delta]^+) 
= H_0(|\Delta|; k) = \rH_0(|\Delta|; k) \oplus k$. 
So Theorem~\ref{cohomology} coincides with Hochster's formula. 
If $0 \ne \ba \in \NN^n$ and 
$\supp(\ba) = F$, 
Hochster's formula states that $[H_\m^{i+1}(k[\Delta])]_{-\ba} = 
H_i(|\Delta|, |\Delta| - \{ p \} ;k)$ for a point $p \in |F|^\circ$. 
Set  $u_F := U_F \cap |\Delta|$. For $p \in |F|^\circ$, $u_F$ is a 
cone neighbourhood of $p$ and $|\Delta| - u_F$ is a deformation retract of 
$|\Delta| - \{p\}$. Hence we have 
\begin{eqnarray*}
[H_\m^{i+1}(M)]_{-\ba} &\cong& 
H_c^i(U_F, j^* k[\Delta]^+) \quad 
\text{(by Theorem~\ref{cohomology2})} \\
&\cong& H_c^i(\, u_F, \, \const_{u_F} \, ) \\ 
&\cong&  H^i(\, |\Delta|, \, |\Delta| - u_F ; \,  k \, ) 
\quad (\text{see \cite[IV. Definition~8.1]{Iver}}) \\
&\cong& 
H^i(\, |\Delta|, \, |\Delta| - \{p\} ; \,  k \, ).
\end{eqnarray*}
So Theorem~\ref{cohomology2} and Hochster's formula also coincide. 
\end{rem}

\section{Relation to Poincar\'e-Verdier Duality}
Since the functor $(-)^+ : \Sq \to \Sh(B)$ is exact, 
it can be extended to the functor $(-)^+ : D^b(\Sq) \to D^b(\Sh(B))$. 
If $\M \in D^b(\Sq)$, we have 
$R\Hom_{k[Q]}(\M, \can) \in D^b_{\Sq}(\MMZn)$ by Lemma~\ref{local dual}.  
So there is a bounded complex $\N$ of squarefree modules such that 
$\N \cong R \Hom_{k[Q]}(\M, \can)$ in $D^b(\MMZn)$. 
We denote $(\N)^+  \in D^b(\Sh(B))$ by $R \Hom_{k[Q]}(\M, \can)^+$. 
Of course, $R \Hom_{k[Q]}(\M, \can)^+$ does not depend on 
the particular choice of $\N$ up to isomorphism in $D^b(\Sh(B))$. 

\medskip
 
For a locally compact topological space $X$ of finite dimension 
(e.g., a locally closed subset of $B$), $\D_X$ denotes a dualizing 
complex of $X$ with the coefficients in $k$ (see \cite[V. \S2]{Iver}). 
In this paper, we frequently use the isomorphism 
$\D_Y \cong j^!\D_X$ for the embedding map $j$ from a locally closed subset 
$Y$ to $X$ (see \cite[V. Theorem~5.6]{Iver}). 
If $X$ is a manifold (with or without boundary), we have the 
orientation sheaf $\Or_X$ of $X$ with the coefficients in $k$. In this case, 
we have $\D_X \cong \Or_X[\dim X]$ (see \cite[V. \S3]{Iver}). 

\begin{lem}\label{Or_B}
With the above notation, we have the following. 

\begin{itemize}
\item[(1)] $\Or_B \cong (\can)^+$.
\item[(2)] Let $J_F$ be the projective object in $\Sq$ associated with a face 
$F \in L$. Then $R\cHom_{\Sh(B)}((J_F)^+, \Or_B) 
\cong \Hom_{k[Q]}(J_F, \can)^+$. 
\item[(3)] If $\M \in D^b(\Sq)$, we have an isomorphism 
$$R\cHom_{\Sh(B)}((\M)^+, \Or_B) \cong 
R\Hom_{k[Q]}(\M, \can)^+$$ in $D^b(\Sh(B)).$ 
\end{itemize}
\end{lem}

\begin{proof}
(1) Let $\const_{B^\circ}$ be the constant sheaf on the relative interior 
$B^\circ$ of $B$. If $j: B^\circ \to B$ is the embedding map, then 
$\Or_B \cong j_! \const_{B^\circ}$ by \cite[VI. Proposition~3.3]{Iver}. 
On the other hand, $(\can)^+ \cong j_! \const_{B^\circ}$ as we have seen in  
Example~\ref{constant}. 

(2) Recall that if $U$ is an open set with the embedding map 
$j: U \to B$ and $\cI$ is an injective object in $\Sh(B)$, then 
$j^* \cI$ $(= j^! \cI)$ is injective in $\Sh(U)$. 
So $\cExt_{\Sh(B)}^i((J_F)^+, \Or_B)$ is the sheaf associated 
to the presheaf which sends an open set $U$ to 
$\Ext^i_{\Sh(U)}(j^* (J_F)^+, j^*\Or_B)$. Note that $j^*\Or_B \cong 
\Or_U$. By the Poincar\'{e}-Verdier duality 
(\cite[V. 2.1]{Iver}), we have 
$$\Ext^i_{\Sh(U)}(j^* (J_F)^+, j^*\Or_B) \cong 
H_c^{n-1-i}(U, j^* (J_F)^+)^\vee,$$
where $(-)^\vee$ means the dual $k$-vector space.
For any open neighbourhood $V$ of $p$, there is an open set $U$ with 
$p \in U \subset V$ such that $U \cap U_F \cong \RR^{n-1}$ or $\RR^{n-1}_+$. 
Then $H_c^i(U, j^* (J_F)^+) \cong H_c^i(U \cap U_F ; k) =0$ for all 
$i \ne n-1$. Thus $\cExt^i_{\Sh(B)}((J_F)^+, \Or_B) = 0$ for all $i \ne 0$. 
Hence we have 
 $R\cHom_{\Sh(B)}((J_F)^+, \Or_B) \cong 
 \cHom_{\Sh(B)}((J_F)^+, \Or_B)$. 

Recall that $(J_F)^+$ is the constant sheaf on $U_F$ and $\Or_B$ 
is the constant sheaf on $B^\circ$. 
For a point $p \in B$, the stalk $\cHom_{\Sh(B)}((J_F)^+, \Or_B)p$ 
at $p$ is nonzero (equivalently, $\cHom_{\Sh(B)}((J_F)^+, \Or_B)p = k$) 
if and only if there is an open neighbourhood $U_p$ of $p$ 
such that $U_p \cap U_F \subset B^\circ$. 
With the same notation as Lemma~\ref{basic} (3), 
the latter condition is equivalent to the condition that 
$\supp(p) \vee F = \RR_{\geq 0}Q$. 
So the assertion follows from Lemma~\ref{basic} (3). 

(3) Let $\P$ be a projective resolution of $\M$ in $\Sq$, 
that is, there is a quasi isomorphism $\P \to \M$ and each $P^i$ is a 
direct sum of copies of $J_F$ for various $F$. 
By (2), we can compute $R \cHom_{Sh(B)}((\M)^+, \Or_B)$ by $(\P)^+$. 
So we have  
\begin{eqnarray*}
R\cHom_{\Sh(B)}((\M)^+, \Or_B) &\cong& 
\cHom_{\Sh(B)}^\bullet((\P)^+, \Or_B) \\
&\cong& \Hom_{k[Q]}^\bullet(\P, \can)^+  \\
&\cong& R\Hom_{k[Q]}(\M, \can)^+.
\end{eqnarray*}
\end{proof}

The normalized $\ZZ^n$-graded dualizing complex of $k[Q]$ is 
a $\ZZ^n$-graded injective resolution of $\can[n]$. 
But, in this paper, we will consider a $\ZZ^n$-graded injective resolution of 
$\can[n-1]$, which is a non-normalized dualizing complex. 
The reason why we use this convention is that 
$k[Q]$ represents the $(n-1)$-dimensional polytope $B$ in our context. 

Let $\DS_{k[Q]}$ be the $Q$-graded part of a minimal $\ZZ^n$-graded 
injective resolution of $K_{k[Q]}[n-1]$. The complex $\DS_{k[Q]}$, 
which is a minimal injective resolution of $\can[n-1]$ in $\Sq$, 
is of the form 
\begin{equation}\label{DC}
\DS_{k[Q]} : 0 \too \omega^{-n+1} \too \omega^{-n+2} \too \cdots \too \omega^1 
\too 0,\end{equation}
$$\omega^i = \bigoplus_{\substack{F \in L \\ \dim F = -i+1}} k[F], $$
and the differential is composed of 
the maps $\varepsilon(F, G) \cdot 
\operatorname{nat} : k[F] \to k[G]$ for all $G \in L$ with  
$\dim G = \dim F -1$, where $\varepsilon$ is the incidence function 
on the cell complex $B = \coprod_{F \in L} |F|^\circ$ and 
$\operatorname{nat} : k[F] \to k[G]$ is the natural surjection. 

 For an order ideal $\Delta \subset L$, set 
$\DS_{k[\Delta]} := \Hom_{k[Q]}(k[\Delta], \DS_{k[Q]})$. 
This is a complex of squarefree $k[\Delta]$-modules with 
$$\omega_{k[\Delta]}^i = \bigoplus_{\substack{F \in \Delta 
\\ \dim F = -i+1}} k[F].$$ 
Note that $\DS_{k[\Delta]}$ is isomorphic to 
a non-normalized $\ZZ^n$-graded dualizing complex of $k[\Delta]$ 
in the derived category of $\ZZ^n$-graded $k[\Delta]$-modules. 

Let $\Sq(\Delta)$ be the full subcategory of $\Sq$ consisting 
of $k[\Delta]$-modules, that is,  
$M \in \Sq(\Delta)$ if and only if $M$ is a squarefree $k[Q]$-module 
whose annihilator $\ann(M)$ contains $I_\Delta$. 
The category $\Sq(\Delta)$ is a thick abelian subcategory of $\Sq$, 
and $\Sq(\Delta)$ has enough injectives, 
and an indecomposable injective object is of the form $k[F]$ for 
some $F \in \Delta$ (c.f. \cite{RWY}), which is also injective in $\Sq$. 
Thus $D^b(\Sq(\Delta)) \cong D^b_{\Sq(\Delta)}(\Sq)  
\cong D^b_{\Sq(\Delta)}(\MMZn)$, and $D^b(\Sq(\Delta))$ can be viewed 
as a full subcategory of $D^b(\MMZn)$. 
 
If $\M \in D^b(\Sq(\Delta))$, we have $R \Hom_{k[Q]}(\M, K_{k[Q]}[n-1]) 
\cong R \Hom_{k[\Delta]}(\M, \DS_{k[\Delta]})$ in $D^b(\MMZn)$ 
by the local duality. 
In particular, $R \Hom_{k[\Delta]}(\M, \DS_{k[\Delta]})$ belongs 
to $D^b_{\Sq(\Delta)}(\MMZn)$, and we can define 
$R \Hom_{k[\Delta]}(\M, \DS_{k[\Delta]})^+ \in D^b(\Sh(B))$. 

If $M \in \Sq(\Delta)$ and $j: |\Delta| \to B$ is the embedding map, 
then $\Supp(M^+) \subset |\Delta|$ and $j_*j^* M^+ \cong M^+$. 
Since $j_* (= j_!) : \Sh(|\Delta|) \to \Sh(B)$ is an exact functor 
in this case, it can be extended to the functor 
$j_* : D^b(\Sh(|\Delta|)) \to D^b(\Sh(B))$.

\begin{thm}\label{local Verdier duality}
With the above notation, for $\M \in D^b(\Sq(\Delta))$, we have  
$$R \cHom_{\Sh(|\Delta|)}(\, j^*(\M)^+, \D_{|\Delta|} \,) \cong
j^* (R\Hom_{k[\Delta]} ( \, \M, \DS_{k[\Delta]})^+)$$
in $D^b(\Sh(|\Delta|))$. 
\end{thm}

\begin{proof}
In $D^b(\Sh(B))$, we have the following isomorphisms. 
\begin{eqnarray*}
&& j_* R \cHom_{\Sh(|\Delta|)}( \, j^*(\M)^+, \D_{|\Delta|} \, ) \\
&\cong&  j_* R \cHom_{\Sh(|\Delta|)}(\, j^*(\M)^+, j^! \D_B \, ) \\
&\cong& R \cHom_{\Sh(B)}( \, j_* j^* (\M)^+, \D_B \,) 
\qquad \text{(by \cite[VII. Theorem~5.2]{Iver})}\\ 
&\cong& R\cHom_{\Sh(B)}( \, (\M)^+, \Or_B[n-1] \, ) \\ 
&\cong& R\Hom_{k[Q]}( \, \M, \can[n-1] \, )^+ \qquad 
\text{(by Lemma~\ref{Or_B} (3))} \\
&\cong& R\Hom_{k[\Delta]} (  \M, \DS_{k[\Delta]}  )^+.
\end{eqnarray*}
Hence $j_* R \cHom_{\Sh(|\Delta|)}(\, j^*(\M)^+, \D_{|\Delta|} \, ) \cong
R\Hom_{k[\Delta]} (\M, \DS_{k[\Delta]})^+$. 
Applying $j^*$ to the both sides of this isomorphism, 
we have the expected isomorphism. 
\end{proof}

\begin{cor}\label{dualizings}
With the above notation, we have
$\D_{|\Delta|} \cong j^*(\DS_{k[\Delta]})^+$. 
\end{cor}

\begin{proof}
$$\D_{|\Delta|} 
\cong R\cHom_{\Sh(|\Delta|)}(\const_{|\Delta|}, \D_{|\Delta|}) 
\cong j^* (R\Hom_{k[\Delta]}(k[\Delta], \DS_{k[\Delta]})^+)
\cong j^* (\DS_{k[\Delta]})^+$$            
\end{proof}

\begin{prop}\label{locally closed}
Let $\Delta, \Sigma \subset L$ be order ideals with 
$\Delta \supset \Sigma$, and $j$ 
the embedding map from $Z := |\Delta| - |\Sigma|$ to $B$. Then 
$$\D_Z \cong
j^* (R\cHom_{k[Q]}(\, I_{\Delta / \Sigma}, \, 
\DS_{k[Q]})^+),$$ 
where $I_{\Delta / \Sigma} := I_\Sigma/I_\Delta$. 
\end{prop}

\begin{proof}
In $D^b(\Sh(B))$, we have the following isomorphisms. 
\begin{eqnarray*}
&& R\Hom_{k[Q]}(\, I_{\Delta / \Sigma}, \, \DS_{k[Q]})^+ \\
&\cong& R\cHom_{\Sh(B)}(\, j_! \const_Z, \, 
\D_{B} \, ) \quad \quad \text{(by Theorem~\ref{local Verdier duality})}\\
&\cong& Rj_* R\cHom_{\Sh(Z)}
(\,  \const_Z, \, j^! \D_{B} \, ) 
\quad \quad \text{(by \cite[VII. Theorem~5.2]{Iver})}\\
&\cong& Rj_* R\cHom_{\Sh(Z)}
(\, \const_Z, \, \D_Z \, )   \\
&\cong& Rj_* \D_Z. 
\end{eqnarray*}
Hence $R\Hom_{k[Q]}(\, I_{\Delta / \Sigma}, \, \DS_{k[Q]})^+ 
\cong Rj_* \D_Z$. Applying $j^*$ to the both sides 
of this isomorphism, we have the expected isomorphism. In fact, since 
the functor $j^*j_* : \Sh(Z) \to \Sh(Z)$ is natural equivalent to the identity 
functor, we have $j^* Rj_* \cong j^* j_* \cong \Id$ as an endofunctor on 
$D^b(\Sh(Z) )$. 
\end{proof}

In our context, the notion of  a Buchsbaum ring is natural and important. 
The original definition of a Buchsbaum ring (see \cite{SV}) is slightly 
complicated, but for $k[\Delta]$, we have a simple criterion. 

\begin{lem}\label{concentrate}
Let $A = \bigoplus_{i \geq 0} A_i$ be a noetherian $\NN$-graded 
commutative ring with the graded maximal ideal 
$\m = \bigoplus_{i > 0} A_i$ (thus $A_0 = k$ is a field). 
Let $M$ be a finitely generated graded $A$-module of dimension $r$. 
If there is some $s \in \ZZ$ such that $[H_\m^i(M)]_t =0$ for all 
$i < r$ and $t \ne s$, then $M$ is a Buchsbaum $A$-module. 
\end{lem}

If $A$ is generated by $A_1$ as a $k$-algebra, the above fact is a special 
case of the well-known result \cite[I. Proposition~3.10]{SV}. 
Even in the general case, this fact was essentially pointed out in 
\cite{SS}.

\begin{proof} Note that $A$ has a graded normalized dualizing complex $\I_A$.
 Set $\N := \tau_{-r} \Hom_A^\bullet(M,\I_A)$. Here, for a complex $\C$, 
$\tau_{-r} \C$ is the truncated complex 
$$\cdots \too 0 \too 
\ker (C^{-r+1} \to C^{-r+2}) \too C^{-r+1} \too C^{-r+2} \too \cdots.$$ 
We have $H^i(\N) = 0$ for all $i \leq -r$, and $H^i(\N)$ is the graded 
$k$-dual of $H_\m^{-i}(M)$ for all $i > -r$ by the local duality. 
So the cohomologies of $\N$ are concentrated in 
the degree $-s$ components. By \cite[II.Theorem~4.1]{SV}, 
it suffices to prove that $\N$ is isomorphic to a complex of $k$-vector 
spaces in the derived category of graded $A$-modules. 
For a graded $A$-module $N$ and an integer $t$, set $N_{\geq t} := 
\bigoplus_{i \geq t} N_i$. Then chain maps $\N_{\geq -s} \to \N$ and 
$\N_{\geq -s} \to \N_{\geq -s}/ \N_{\geq -s+1}$ are quasi-isomorphisms. 
Thus, in the derived category,  $\N$ is isomorphic to 
$\N_{\geq -s}/ \N_{\geq -s+1}$, which is a complex of $k$-vector spaces.  
\end{proof}

\begin{cor}[c.f. {\cite[Corollary~3.8]{Y2}}]\label{Bbm module}
Let $M$ be a squarefree $k[Q]$-module of dimension $r$. 
Then the following are equivalent. 
\begin{itemize}
\item[(a)] $M$ is a Buchsbaum module. 
\item[(b)] $\dim_k H_\m^i(M) < \infty$ for all $i \ne r$. 
\item[(c)] $[H_\m^i(M)]_\ba = 0$, if $i \ne r$ and  $\ba \ne 0$. 
\end{itemize}
\end{cor}

\begin{proof}
The implication (a) $\Rightarrow$ (b) is a basic property of 
Buchsbaum modules. 
The $\ZZ^n$-graded $k$-dual of $H_\m^i(M)$ is the squarefree 
module $\Ext^{n-i}_{k[Q]}(M, \can)$. Hence 
$\dim_k H_\m^i(M) < \infty$ if and only if $H_\m^i(M) 
= [H_\m^i(M)]_0$. So we have (b) $\Leftrightarrow$ (c). 
The implication (c) $\Rightarrow$ (a) follows from Lemma~\ref{concentrate}. 
\end{proof}

\begin{cor}\label{Bbm}
Let $\Delta \subset L$ be an order ideal with $d = \dim |\Delta|$ 
(so $\dim k[\Delta] = d+1$). The following are equivalent. 
\begin{itemize}
\item[(a)] $k[\Delta]$ is a Buchsbaum ring. 
\item[(b)] $\cH^i(\D_{|\Delta|}) = 0$ for all $i \ne - d$. 
\item[(c)] $H_i(|\Delta|, |\Delta| - \{p\}; \, k) = 0$ 
for all $i < d$ and all $p \in |\Delta|$. 
\end{itemize}
In particular, the Buchsbaum property of $k[\Delta]$ is a topological 
property of $|\Delta|$ (i.e., depends only on the topology of 
$|\Delta|$ and $\chara(k)$). 
\end{cor}

\begin{proof}
We have $\cH^i(\D_{|\Delta|}) \cong 
j^*(H^i(\DS_{k[\Delta]})^+) \cong j^*(\Ext_{k[Q]}^i(k[\Delta], 
\DS_{k[Q]})^+)$ by Corollary~\ref{dualizings},  
where $j : |\Delta| \to B$ is the embedding map. 
Thus $\cH^i(\D_{|\Delta|}) = 0$ if and only if 
$\dim_k \Ext_{k[Q]}^i(k[\Delta], \DS_{k[Q]}) 
= \dim_k H_\m^{-i+1}(k[\Delta])  < \infty$. 
So (b) is equivalent to (a). The equivalence (b) $\Leftrightarrow$ (c) 
must be obvious for algebraic topologists. But the equivalence 
 (c) $\Leftrightarrow$ (a) ($\Leftrightarrow$ (b))
also follows from Theorem~\ref{cohomology2}. In fact, if 
$0 \ne \ba \in \NN^n$, 
we have $[H^{i+1}_\m(k[\Delta])]_{-\ba} \cong 
H_i(|\Delta|, |\Delta| - \{ p \}; k)$ 
for $p \in  | \supp(\ba)|^\circ$, 
as we have seen in Remark~\ref{Hoch}. 
\end{proof}

\begin{rem}\label{Bbm relative}
The implication (b) $\Rightarrow$ (a) of Corollary~\ref{Bbm} does not hold for 
a locally closed subset $Z:= |\Delta| - |\Sigma|$ and its 
squarefree module $I_{\Delta / \Sigma} := I_\Sigma/I_\Delta$, 
while we have Proposition~\ref{locally closed}. 
Since $$R\Hom_{k[Q]}(I_{\Delta / \Sigma}, \DS_{k[Q]})^+ \cong 
Rj_* \D_Z$$ by the proof of Proposition~\ref{locally closed}, 
$I_{\Delta / \Sigma}$ is a Buchsbaum module of dimension $d+1$ if and only if 
$R^ij_* \D_Z = 0$ for all $i \ne -d$. 
These conditions are stronger than the condition that 
$\cH^i(\D_Z) = 0$ for all $i \ne -d$. 

For example,  consider a polynomial ring $k[x,y,z]$, and simplicial  
complexes $\Delta = 2^{\{x,y,z\}}$ and $\Sigma = \{ \{ x \}, \emptyset \}$. 
Then $Z$ is a manifold with boundary (in fact, $Z \cong \RR_+^2$), 
and $\cH^i(\D_Z) = 0$ for all $i \ne -2$. But $I_{\Delta / \Sigma} = 
(y,z)$ is not a Buchsbaum module. In this case, $R^ij_*\D_Z \ne 0$ 
for $i = -1, -2$. 

Since $\Supp ( R^ij_* \D_Z)  \subset \bar{Z} = |\Delta|$, 
it suffices to check $R^ih_* \D_Z$ to see the vanishing of $R^ij_* \D_Z$, 
where $h: Z \to |\Delta|$ is the embedding map. That is, 
$I_{\Delta / \Sigma}$ is a Buchsbaum module of dimension $d+1$ 
if and only if $R^ih_* \D_Z = 0$ for all $i \ne -d$. 
Hence the Buchsbaum property of $I_{\Delta / \Sigma}$ is a topological 
property of the pair $(|\Delta|, |\Sigma|)$. 
\end{rem}

If $|\Delta|$ is a manifold (with or without boundary) of dimension $d$, 
then we have  $\D_{|\Delta|} \cong \Or_{|\Delta|}[d]$ and 
$k[\Delta]$ is a Buchsbaum ring of dimension $d+1$. 
By Corollary~\ref{dualizings}, we have 
$j^* (K_{k[\Delta]})^+ \cong \Or_{|\Delta|}$, where 
$K_{k[\Delta]} := \Ext_{k[Q]}^{n-d-1}(k[\Delta], \can)$ is the 
canonical module of $k[\Delta]$. 

Let $(A, \m)$ is a Buchsbaum local ring of dimension $d+1$ 
admitting a canonical modules $K_A$. 
Then \cite[II. Theorem~4.9]{SV} states that  
$$H_\m^i(K_A) \cong \Hom_A(H_\m^{d-i+2}(A), \, E(A/\m) ) 
\quad \text{for all $2 \leq i \leq d$,}$$
where $E(A/\m)$ is the injective hull of $A/\m$. 
We will see that this duality corresponds to the Poincar\'e duality 
in our context. 

Assume that $k[\Delta]$ is a Buchsbaum ring of dimension $d+1$ 
(thus $\dim |\Delta| = d$). Then  we have 
\begin{equation}\label{Poincare}
[H_\m^i(\, (K_{k[\Delta]})_{> 0} \, )]_0 \cong 
[H_\m^{d-i+2}(\, k[\Delta]_{>0} \, )^\vee]_0 
\quad \text{for all $1 \leq i \leq d+1$.}
\end{equation}
(When $2 \leq i \leq d$, this is just a $\ZZ^n$-graded version 
of \cite[II. Theorem~4.9]{SV}. We leave the case when $i = 1, d+1$ 
for the reader as an easy exercise.)
By Theorem~\ref{cohomology}, 
$$[\, H_\m^i(\, (K_{k[\Delta]})_{> 0} \,) \, ]_0 \cong 
H^{i-1}(\, |\Delta|, \, (K_{k[\Delta]})^+) \cong 
H^{i-1}(\, |\Delta|, \, \Or_{|\Delta|} \, )$$
and  
$$[\, H_\m^{d-i+2}(\, k[\Delta]_{>0} \, ) \, ]_0 \cong 
H^{d-i+1}(\, |\Delta|, \, \const_{|\Delta|} \, ) \cong 
H^{d-i+1}(\, |\Delta| \, ; \, k \, )$$ 
for all $1 \leq i \leq d+1$. So \eqref{Poincare}  
also follows from the Poincar\'e duality 
$$H^i(\, |\Delta|, \, \Or_{|\Delta|} \, ) \cong 
H^j(\, |\Delta| \, ; \, k \, )^\vee \ 
\text{for all $i, j$ with $i+j = d$.}$$
Note that $|\Delta|$ is an orientable manifold (i.e., a manifold with 
$\const_{|\Delta|} \cong \Or_{|\Delta|}$) if and only if $k[\Delta]$ 
is a Buchsbaum ring with $(K_{k[\Delta]})_{>0} \cong k[\Delta]_{>0}$. 
In this case, \eqref{Poincare} corresponds to the most 
familiar form of the Poincar\'e duality.  
We also remark that if $|\Delta|$ is an orientable manifold 
of dimension $d$ then $\dim_k [H_\m^{d+1}(k[\Delta])]_0$ equals the number of 
the connected components of $|\Delta|$. When $|\Delta|$ is a 
connected manifold, $|\Delta|$ is orientable if and only if 
$\dim_k [H_\m^{d+1}(k[\Delta])]_0 = 1$. 
In this case, $K_{k[\Delta]} \cong k[\Delta]$.

\medskip

Let $\Sq_+(\Delta)$ be the full subcategory of $\Sq$ consisting of 
squarefree $k[\Delta]$-modules $M$ with $M_0=0$. For a while, let 
$\M$ be an object of $D^b(\Sq_+(\Delta))$. 

For $\M \in D^b(\Sq_+(\Delta))$, 
by the local duality and Theorem~\ref{cohomology}, we have 
$$[\Ext^{-i}_{k[\Delta]}(\M, \DS_{k[\Delta]})^\vee]_0 
\cong [ R^{i+1}   \Gamma_\m(\M) ]_0 
\cong R^i \, \Gamma (B, (\M)^+ ) 
\cong R^i \, \Gamma ( \, |\Delta|, \, j^*(\M)^+).$$
On the other hand, we have 
$\Ext^{-i}_{\Sh(|\Delta|)}(\, j^*(\M)^+, \, \D_{|\Delta|} \, )^\vee  
\cong R^i \, \Gamma( \, |\Delta|, \, j^*(\M)^+ ) $
by the  Poincar\'e-Verdier duality (\cite[V, 2.1]{Iver}). Thus 
\begin{equation}\label{global duality}
\Ext^i_{\Sh(|\Delta|)}(\, j^*(\M)^+, \, \D_{|\Delta|} \, )  
\cong [\Ext^i_{k[\Delta]}(\, \M, \, \DS_{k[\Delta]} \, )]_0.
\end{equation}

We can give another proof of \eqref{global duality}. 
Let $\P \to \M$ be a projective resolution in $\Sq$. 
Since $\M \in D^b(\Sq_+(\Delta))$, we may assume that 
each component of $\P$ is a direct sum of copies of 
$J_F$ for various $\{ 0 \} \ne F \in L$. If $F \ne \{ 0 \}$, 
then $\Supp( (J_F)^+) = U_F \cong \RR^{n-1}$ or $\RR^{n-1}_+$ and 
$\Ext^i_{\Sh(B)}( (J_F)^+, \, \Or_B  ) = H^{n-1-i}(B, (J_F)^+) = 
H_c^{n-1-i}(U_F; k) = 0$ for all $i \ne 0$. So we can compute 
$\Ext^i_{\Sh(B)}(\, (\M)^+, \, \Or_B \,)$
using $(\P)^+$, and we have the following. 
\begin{eqnarray*}
& & \Ext^i_{\Sh(|\Delta|)}(\, j^*(\M)^+, \, \D_{|\Delta|} \, ) \\
&\cong& \Ext^i_{\Sh(B)}(\, (\M)^+, \, \D_B \,) 
\qquad  \text{(by \cite[VII. Theorem~3.1]{Iver})} \\
&\cong& \Ext^i_{\Sh(B)}(\, (\M)^+, \, \Or_B[n-1] \, ) \\
&\cong& H^i(\, \Hom_{\Sh(B)}(\, (\P)^+, \, (K_{k[Q]})^+[n-1] \, ) \, ) \\
&\cong& H^i(\, [ \, \Hom_{k[Q]}(\, \P, \, K_{k[Q]}[n-1] \,) \, ]_0 \,) 
\qquad  \text{(since $\P \in D^b(\Sq_+)$)} \\
&\cong& [\, \Ext^i_{k[Q]}(\, \M, \, K_{k[Q]}[n-1] \, ) \, ]_0 \\
&\cong& [\, \Ext^i_{k[\Delta]}(\, \M, \, \DS_{k[\Delta]} \, ) \, ]_0.
\end{eqnarray*}  

Finally, we study the Cohen-Macaulay property of $k[\Delta]$ and 
$I_{\Delta/\Sigma}$. 
If $\dim k[\Delta] \leq 1$, then $k[\Delta]$ is always Cohen-Macaulay. 
So we may assume that $\dim k[\Delta] \geq 2$. 
The same thing is true for $I_{\Delta/\Sigma}$. 
When $k[Q]$ is a polynomial ring, the next result is a well-known 
theorem of Munkres.

\begin{thm}[c.f. \cite{Mil, Y3}]
Let $\Delta \subset L$ be an order ideal with 
$d := \dim |\Delta| \geq 1$ (i.e., $\dim k[\Delta] = d+1 \geq 2$). 
Then the following are equivalent.

\begin{itemize}
\item[(a)] $k[\Delta]$ is a Cohen-Macaulay ring of dimension $d+1$,
\item[(b)] $\rH_i(|\Delta| ; k) 
= H_i(|\Delta|, |\Delta|-\{ p \} ; k) = 0$ for all 
$i < d$ and all $p \in |\Delta|$,
\item[(c)] $\cH^i (\D_{|\Delta|}) = 0$ for all $i \ne -d$, 
$H^i \Gamma( |\Delta|, \D_{|\Delta|}) = 0$ for all $i \ne -d, \, 0$, 
and $H^0 \Gamma(|\Delta|, \D_{|\Delta|}) \cong k$. 
\end{itemize}
In particular, the Cohen-Macaulay property of $k[\Delta]$ 
is a topological property of $|\Delta|$. 
\end{thm}

\begin{proof}
The equivalence between (a) and (b) has been proved in \cite{Mil, Y3}. 
Recall that $H_i(|\Delta|, |\Delta|-\{ p \} ; k) = 0$ 
for all $i < d$ and all $p \in |\Delta|$ if and only if 
$\cH^i (\D_{|\Delta|}) = 0$ for all $i \ne -d$. Since 
$H^{-i} \Gamma( |\Delta|, \D_{|\Delta|}) 
\cong H^i(|\Delta| ; k)^\vee$, (b) and (c) are equivalent. 
\end{proof}

\begin{prop}
Let $\Delta, \Sigma \subset L$ be order ideals with 
$\Delta \supset \Sigma \ne \emptyset$, and $h$ the embedding map 
from $Z := |\Delta| - |\Sigma|$ to $|\Delta|$.  

\begin{itemize}
\item[(a)] $I_{\Delta/\Sigma}$ is a Cohen-Macaulay module of dimension $d+1$, 
\item[(b)] $R^ih_*\D_Z = H^i \Gamma(Z, \D_Z) = 0$ for all $i \ne -d$. 
\end{itemize}
In particular, the Cohen-Macaulay property of $I_{\Delta/\Sigma}$ 
is a topological property of the pair $(|\Delta|, |\Sigma|)$. 
\end{prop}

\begin{proof}
$I_{\Delta/\Sigma}$ is Cohen-Macaulay if and only if it is Buchsbaum 
and $[H^i_\m(I_{\Delta/\Sigma})]_0 = 0$ for all $i \ne d+1$. 
As we have seen in Remark~\ref{Bbm relative}, $I_{\Delta/\Sigma}$ is Buchsbaum 
if and only if $R^ih_*\D_Z = 0$ for all $i \ne -d$. 
Since $I_{\Delta/\Sigma} \in \Sq_+$, we have 
$[H^{i+1}_\m(I_{\Delta/\Sigma})]_0 
\cong H^i(B, (I_{\Delta/\Sigma})^+) \cong H^i_c(Z;k) \cong 
H^{-i} \Gamma(Z, \D_Z)^\vee$ for all $i$. So we are done.  
\end{proof}

\end{document}